\documentclass[10pt,twoside]{uz_kgu_my}
\usepackage[cp1251]{inputenc}
\usepackage[english,russian]{babel}

\begin{document}
\renewcommand{\figurename}{Figure}
\renewcommand{\refname}{\small {References}}
\renewcommand{\proofname}{ {\hskip\parindent \bf Proof. }}
\UDK{517.938} \ArticleNAME{ON MEASURE INVARIANCE\\ FOR A $2$-VALUED
TRANSFORMATION} \ArticleAUTHOR{P.I. TROSHIN} \ArticleHEAD{On measure
invariance for a $2$-valued transformation} \ArticleAUTHORHEAD{P.I.
Troshin} \makeabstitle
\vspace{\baselineskip}
\renewcommand{\abstractname}{Summary}
\begin{abstract}
\vspace{3pt} We consider a family $S=S(a)$ of $2$-valued
transformations of special form on the segment $[0,1]$ with measure
$\mu={\displaystyle\int} p(x)\,d\lambda$, which is absolutely
continuous with respect to the Lebesgue measure $\lambda$. We endow
$S$ with a set of weight functions
$\alpha=\{\alpha_1(x),\alpha_2(x)\}$ and find a criterion of measure
invariance under the transformation. This criterion relates the
three parameters $a$, $p$, $\alpha$ to each other.

 \vspace{3pt}
\textbf{Key words:} multivalued dynamical system, invariant measure,
$\beta$-expansion, Bernoulli convolution
\end{abstract}
\renewcommand{\abstractname}{Аннотация}
\begin{abstract}
Рассматривается семейство двузначных трансформаций $S=S(a)$
специального вида на отрезке $[0,1]$ с мерой
$\mu={\displaystyle\int} p(x)\,d\lambda$, абсолютно непрерывной
относительно меры Лебега $\lambda$. Трансформация $S$ оснащается
набором весовых функций $\alpha=\{\alpha_1(x),\alpha_2(x)\}$.
Находится критерий инвариантности меры под действием заданной
оснащенной трансформации. Этот критерий явным образом связывает три
параметра: $a$, $p$ и $\alpha$.

\vspace{3pt} \textbf{Ключевые слова:} многозначная динамическая
система, инвариантная мера, $\beta$-разложение
\end{abstract}
\vspace{\baselineskip}\hrule
\vspace{\baselineskip}

\section{Introduction. Dynamical system connected to arithmetic representation}\label{sec2.1}

Connections between the ergodic theory and the metric number theory
are well known. One of these connections are arithmetic
representations arising in special symbolic realization of dynamical
systems.

Let $\beta\in(1,2]$. Any infinite sequence $\sigma_1\sigma_2\ldots$
of zeros and ones is called {\it$\beta$-expansion} of number
$x\in[0,1]$ (see~\cite{renyi1957,parry1960}), provided that
$$x=\sum_{k=1}^\infty \sigma_k \beta^{-k}.$$ It is clear that with
$\beta=2$ we obtain the usual binary representation of number $x$.

Every number $x\in[0,1]$ has at least one $\beta$-expansion which is
called  {\it canonical}, or {\it<<greedy expansion>>}:
$\sigma_k=[\beta T^{k-1}x]$, $k\ge1$, where $Tx=\{\beta x\}$ ($[y]$
and $\{y\}$ --- whole and fractional parts of number $y\in\mathbb
R$) (see~\cite{sidorov2003b}).

If $\beta\in\left(1,\frac{1+\sqrt{5}}{2}\right)$, then $x\in[0,1]$
has a continuum of different $\beta$-expansions~\cite{erdosh1990},
if $\beta\in\left[\frac{1+\sqrt{5}}{2},2\right)$, the same is true
for almost every $x\in[0,1]$ ~\cite{sidorov2003}. On the other hand,
for all $m\in \mathbb N$ there exists a base
$\beta\in\left(\frac{1+\sqrt{5}}{2},2\right)$ and a number
$x\in[0,1]$ which has exactly $m$ different
$\beta$-expansions~\cite{sidorov2009}.

\begin{figure}[!htb]
\centering\includegraphics[height=3.5cm]{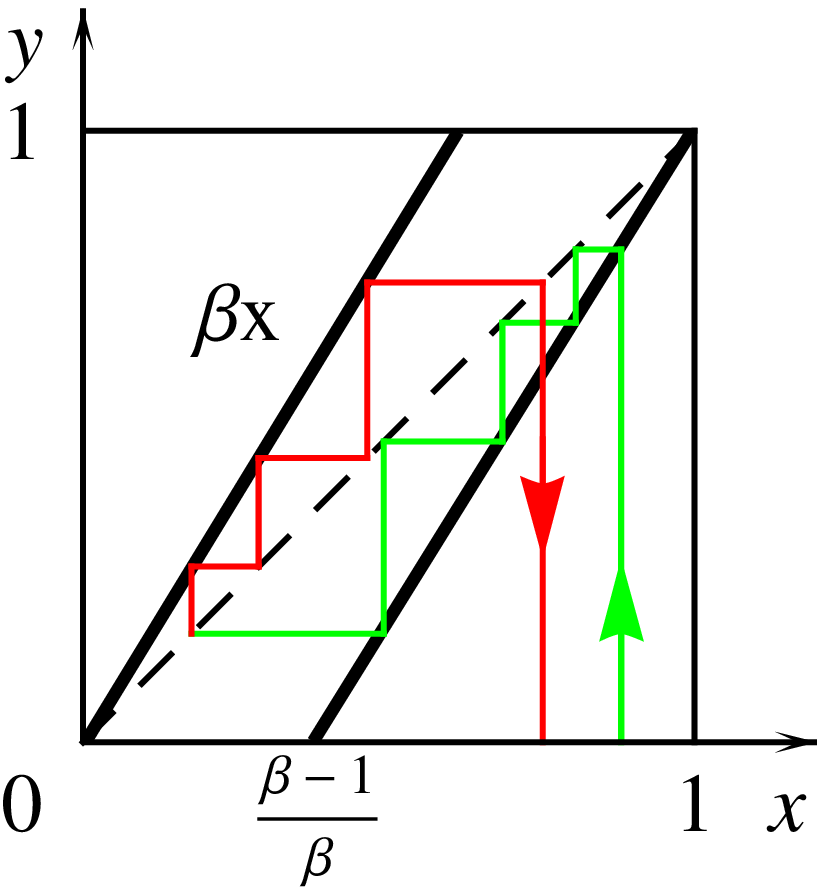}\qquad\qquad
\includegraphics[height=3.5cm]{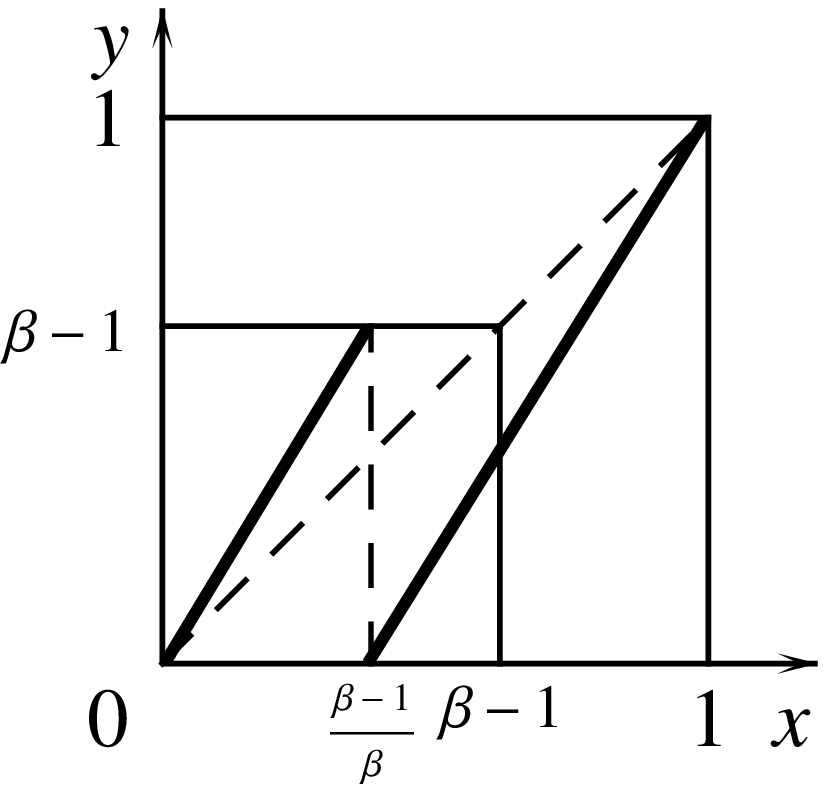} \caption{Scheme of transformations $S$ (on the left) and $S_2$ (on the right)}\label{ergodpic}
\end{figure}

As it is shown in~\cite{KB2009}, to find an arbitrary
$\beta$-expansion of number  $x_1\in[0,1]$, it is  necessary and
sufficient to follow the one of the possible orbit of the point
$x_1$ under multi-valued transformation  $S$ (see
fig.~\ref{ergodpic})
$$\begin{cases}h_0(x)=\beta x,& x\in\left[0,\frac{1}{\beta}\right];\\[0.2cm]
h_1(x)=\beta
x-\beta+1,&x\in\left[1-\frac{1}{\beta},1\right].\end{cases}$$ There
is a correspondence between every orbit $x_1, x_2,x_3\ldots$ and
$\beta$-expansion $\sigma_1\sigma_2\sigma_3\ldots$ of number $x_1$
which can be found by the rule:
$$\sigma_{k}=0,\mbox{ если }
x_{k+1}=h_0(x_{k});$$
$$\sigma_{k}=1,\mbox{ если } x_{k+1}=h_1(x_{k}).$$
If $x_{k}\in[1-\frac{1}{\beta},\frac{1}{\beta}]$, then we can choose
the mapping $h_0$ as well as $h_1$ to construct $x_{k+1}$. Thus we
obtain every possible $\beta$-expansions of number $x_1$.

To find a canonical $\beta$-expansion it is necessary and sufficient
to follow in the same way as before a single-valued orbit under
$1$-transformation $S_2$ (see fig.~\ref{ergodpic})
$$S_2(x)=\begin{cases}\beta x,&
x\in\left[0,1-\frac{1}{\beta}\right);\\[0.2cm]
\beta x-\beta+1,&x\in\left[1-\frac{1}{\beta},1\right].\end{cases}$$

By investigating the orbits of points under transformation $S_2$, we
come to the conclusion that orbits of all the points except $x=1$
are <<captured>>
 by the segment $[0,\beta-1]$. Dynamical system $([0,\beta-1];
S_2|_{[0,\beta-1]})$ was considered in classical papers on
$\beta$-expansions~\cite{renyi1957, gelfond1959, parry1960}, in
which they found invariant measure equivalent to the Lebesgue
measure, and also in~\cite{barnsley2005, KB2009}, where they
calculated top addresses for the iterated function system
$\left\{[0,1];
\phi_1(x)=\frac{1}{\beta}x+1-\frac{1}{\beta},\phi_2(x)=\frac{1}{\beta}x,\quad 1<\beta<2 \right\}$. 

Transformation $S$ is also tightly connected to the problem of
finding every parameter $\beta\in (1,2)$ which provides
nonsingularity of {\it Erd\H{o}s measure} on the segment $[0,1]$,
--- а measure corresponding to a distribution of  random variable
$\sum_{k=1}^\infty \sigma_k \beta^{-k}$ where coefficients
$\sigma_k\in\{0,1\}$ are independently chosen with probability $1/2$
(this is called  {\it Bernoulli convolution problem} and it has not
been solved yet) (see~\cite{erdosh1939}, fine survey on the topic
can be found in~\cite{convolutions2000}).

\begin{figure}[!tbh]
\centering
\begin{minipage}[c]{3.5cm}
\centering\includegraphics[height=3cm,viewport=73 70 295
295,clip]{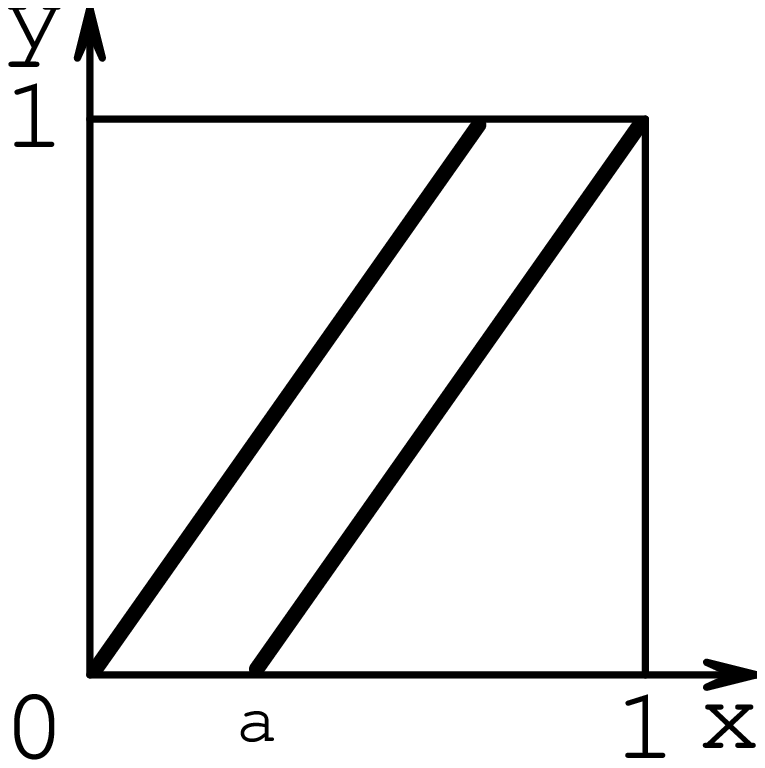}
\end{minipage}%
\begin{minipage}[c]{0.05\textwidth}\centering$=$
\end{minipage}%
\begin{minipage}[c]{3.5cm}
\centering\includegraphics[height=3cm,viewport=25 70 295
295,clip]{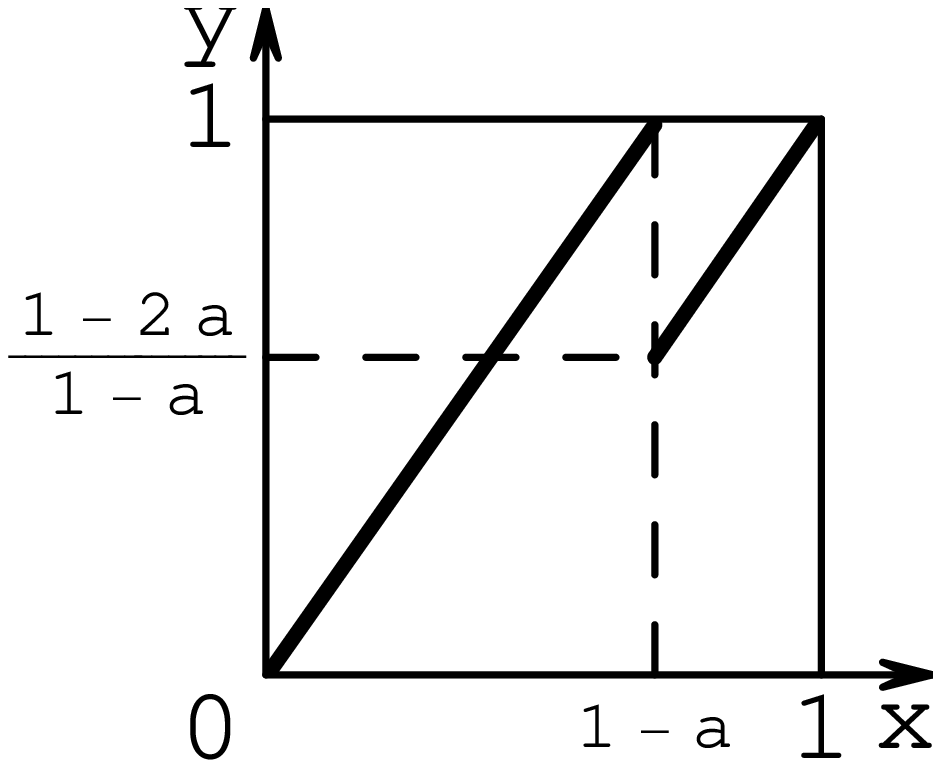}
\end{minipage}%
\begin{minipage}[c]{0.05\textwidth}
\centering$\bigcup$
\end{minipage}%
\begin{minipage}[c]{3.5cm}
\centering\includegraphics[height=3cm,viewport=35 70 295
295,clip]{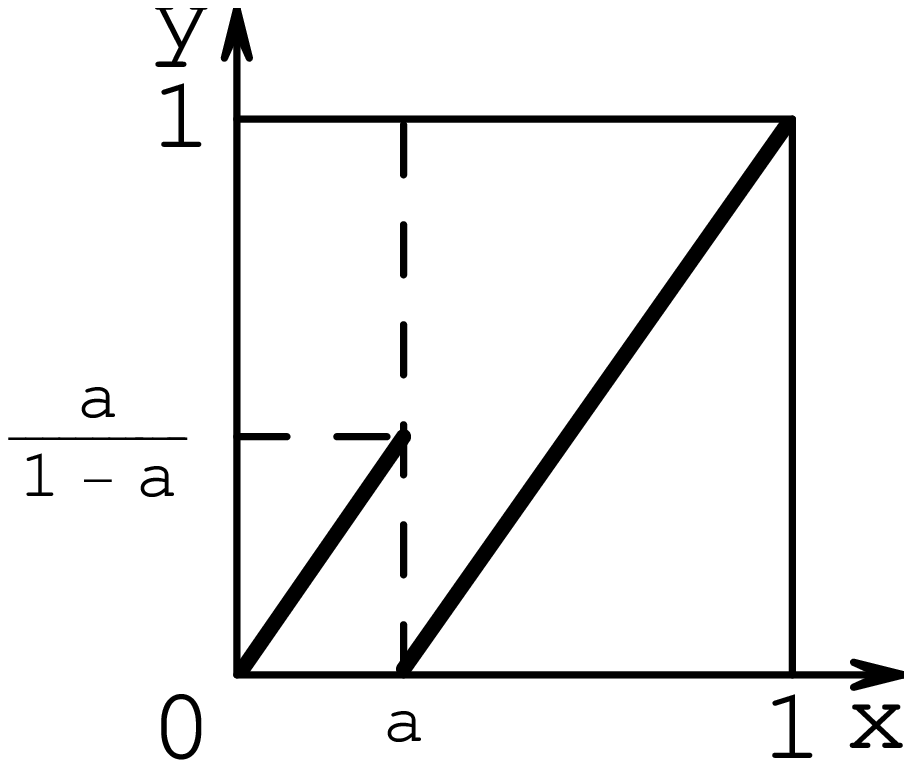}
\end{minipage}%
\caption{Scheme of $2$-transformation $S=S_1\cup
S_2$}\label{picex222}
\end{figure}

Given $a\in \left.\left(0,\frac{1}{2}\right.\right]$, we consider
$2$-transformation $S=S_1\cup S_2$ on the segment $[0,1]$ (see fig.~\ref{picex222}) where\\[2mm]
\begin{tabular}{ll}
$S_1(x)=\begin{cases}\frac{1}{1-a}x,& x\in[0,1-a);\cr
\frac{1}{1-a}x-\frac{a}{1-a},&x\in[1-a,1],\cr\end{cases}$& $\quad
S_2(x)=\begin{cases}\frac{1}{1-a}x,& x\in[0,a);\cr
\frac{1}{1-a}x-\frac{a}{1-a},&x\in[a,1].\cr\end{cases}$
\end{tabular}

\vskip0.2cm\noindent Note that in the discussion above
 $\beta=\frac{1}{1-a}$.


Let $\lambda$ be the Lebesgue measure on $[0,1]$, $\mathfrak B$
--- $\sigma$-field of Borel subsets of  $[0,1]$. Let also $\mu(B)=\int_B p(x)\,d\lambda$ be
a measure absolutely continuous with respect to the Lebesgue
measure, $p(x)\in L^1([0,1], \mathfrak B,\lambda)$ and $p(x)\ge0$.
We endow $2$-transformation $S$ with a set of weight functions
$$\alpha=\{\alpha_1(x),\alpha_2(x)\},\quad \alpha_1(x),\alpha_2(x)\in
L^1([0,1], \mathfrak B,\lambda)$$ provided that
$\alpha_1(x)+\alpha_2(x)\equiv1$ and $\alpha_1(x),\alpha_2(x)\ge0$.
Henceforth, we  consider  $2$-valued dynamical system
$([0,1],\mathfrak {B},\mu, \{S,\alpha\})$ with the finite measure
$\mu$.

According to~\cite{troshin2009}, we obtain a new measure $\mu_S$ on
$\mathfrak B$ by the next formula:
$$\mu_S(B)=\int\limits_{S_1^{-1}(B)}\alpha_1(x)p(x)\, d\lambda+\int\limits_{S_2^{-1}(B)}\alpha_2(x)p(x)\, d\lambda.$$

Denote $\alpha_1(x)p(x)=A_1(x)$, $\alpha_2(x)p(x)=A_2(x)$. Then
$A_1(x)+A_2(x)=p(x)$ and
$$\mu_S(B)=\int\limits_{S_1^{-1}(B)}A_1(x)\,
d\lambda+\int\limits_{S_2^{-1}(B)}A_2(x)\, d\lambda.$$

There are three independent parameters in the studied construction:
density function  $p(x)$, number $a$ (parameter of transformation
$S$) and an equipment $\alpha=\{\alpha_1(x),\alpha_2(x)\}$
($\mu=\mu(p)$ and $S=S(a,\alpha)$).
When searching for equipped transformation with given invariant
measure $\mu$ or searching measure which is preserved by given
transformation $S$,
--- we have a certain relation between the parameters to fulfill the equation $\mu_S=\mu$.

This relation is investigated in Lemma~\ref{lemma1bis} and
Theorem~\ref{main}. In the Corollary~\ref{lebesguemeasure} we
particularly discuss the Lebesgue measure($p\equiv 1$). The
construction defined in Theorem~\ref{interesting} gives an example
of an invariant measure with non-constant density. Two more examples
consistent with classical results in ergodic theory are given in
Corollaries~\ref{diadic} and~\ref{invmeasurecorollary}.


\section{Main results. Criterion of measure invariance}\label{sec2.2}
Fix three parameters: $a\in
\left.\left(0,\frac{1}{2}\right.\right]$,
$\alpha=\{\alpha_1(x),\alpha_2(x)\}$ and $p(x)$. Then the condition
of measure invariance $\mu_S=\mu$ is characterized by the following

\begin{Lemma}\label{specialmeasure}
$\mu_S=\mu$ if and only if for almost every $x\in [0,1]$ (with
respect to  $\lambda$)
\begin{multline}\label{lemma1bis}A_1((1-a)x)+\chi_{\left[\frac{1-2a}{1-a},1\right]}(x)A_1((1-a)x+a)+A_2((1-a)x+a)+\\+\chi_{\left[\left.0,\frac{a}{1-a}\right.\right)}(x)A_2((1-a)x)=\frac{p(x)}{1-a}.\end{multline}
\end{Lemma}

\begin{proof} Let $C_1=\left[\frac{1-2a}{1-a},1\right]$, $C_2=\left[\left.0,
\frac{a}{1-a}\right.\right)$. Then by changing  variables in the
Lebesgue integral we obtain that for every $B\in\mathfrak B$
\begin{multline}\label{troshinnew}
\mu_S(B)=\int\limits_{(1-a)B}A_1(x)\,d\lambda+\int\limits_{(1-a)(B\cap
C_1)+a}A_1(x)\,d\lambda+\int\limits_{(1-a)(B\cap
C_2)}A_2(x)\,d\lambda+\\+\int\limits_{(1-a)B+a}A_2(x)\,d\lambda=\\=(1-a)\left(\int\limits_{B}A_1((1-a)x)\,d\lambda+\int\limits_{B}\chi_{C_1}
(x)A_1((1-a)x+a)\,d\lambda+\right.\\\left.+\int\limits_{B}\chi_{C_2}(x)A_2((1-a)x)\,d\lambda+\int\limits_{B}A_2((1-a)x+a)\,d\lambda\right).
\end{multline}

If $\mu_S=\mu$, then considering arbitrariness of $B\in\mathfrak B$
formula~\eqref{troshinnew} implies~\eqref{lemma1bis}. And
conversely, by substituting the equality~\eqref{lemma1bis}
into~\eqref{troshinnew}, we get $\mu_S=\mu$.
\end{proof}


Let henceforward $$\frac{1}{n+1}<a\le\frac{1}{n}\qquad (n\in \mathbb
N, n\ge2).$$

\begin{Theorem}\label{main}
$\mu_S=\mu$ if and only if the following conditions hold true:
\begin{equation}\label{eqth1}\sum\limits_{k=-1}^{n-1}p(x_0+ka)=\frac{1}{1-a}\sum\limits_{k=-1}^{n-2}p\!\left(\frac{x_0+ka}{1-a}\right), \quad
x_0\in[a,1-(n-1)a);\end{equation}

\begin{equation}\label{eqth2}\sum\limits_{k=-1}^{n-2}p(x_1+ka)=\frac{1}{1-a}\sum\limits_{k=-1}^{n-3}p\!\left(\frac{x_1+ka}{1-a}\right),\quad
x_1\in[1-(n-1)a,2a);\end{equation}

\begin{equation}\label{eqth3}\alpha_1(x+ma)p(x+ma)=\sum\limits_{k=-1}^{m}p(x+ka)-\frac{1}{1-a}\sum\limits_{k=-1}^{m-1}p\!\left(\frac{x+ka}{1-a}\right),\end{equation}
where $x+ma\in [(m+1)a,(m+2)a)$ for  $m=\overline{0,n-3}$,
$x+(n-2)a\in[(n-1)a,1-a)$.

There is no restriction on function $\alpha_1(x)$  on the intervals
$[0,a)$ and $[1-a,1]$.
\end{Theorem}
\begin{proof} Let us consider two cases.

1. First let  $0<a\le \frac{1}{3}$. Then
$\frac{a}{1-a}\le\frac{1-2a}{1-a}$ and formula~\eqref{lemma1bis} can
be written in the following form:
\begin{equation}\label{th1bis}\frac{p(x)}{1-a}=\begin{cases}p((1-a)x)+A_2((1-a)x+a),&
x\in\left[\left.0,\frac{a}{1-a}\right.\right);\\[0.2cm]
A_1((1-a)x)+A_2((1-a)x+a),&x\in\left[\left.\frac{a}{1-a},\frac{1-2a}{1-a}\right.\right);\\[0.2cm]
p((1-a)x+a)+A_1((1-a)x),& x\in\left[\frac{1-2a}{1-a},1\right].
\end{cases}\end{equation}

Consider the first equation of the system:
$\frac{p(x)}{1-a}=p((1-a)x)+A_2((1-a)x+a)$,
$x\in\left[\left.0,\frac{a}{1-a}\right.\right)$. Make a change
$y=(1-a)x+a$, then, for  $y\in[a,2a)$,
$$\frac{1}{1-a}\;p\!\left(\frac{y-a}{1-a}\right)=p(y-a)+A_2(y)=p(y-a)+p(y)-A_1(y).$$

Consider the second equation of the system:
$\frac{p(x)}{1-a}=A_1((1-a)x)+A_2((1-a)x+a)$,
$x\in\left[\left.\frac{a}{1-a},\frac{1-2a}{1-a}\right.\right)$. Make
a change $y=(1-a)x$, then, for $y\in[a,1-2a)$,
$$\frac{1}{1-a}\;p\!\left(\frac{y}{1-a}\right)=A_1(y)+A_2(y+a)=p(y+a)+A_1(y)-A_1(y+a).$$

Consider the third equation of the system:
$\frac{p(x)}{1-a}=p((1-a)x+a)+A_1((1-a)x)$,
$x\in\left[\frac{1-2a}{1-a},1\right]$. Make a change $y=(1-a)x$,
then, for $y\in[1-2a,1-a]$,
$$\frac{1}{1-a}\;p\!\left(\frac{y}{1-a}\right)=p(y+a)+A_1(y).$$

Thus we obtain:
\begin{equation}\label{lemma2bis}\begin{cases}A_1(y)=p(y-a)+p(y)-\frac{1}{1-a}\;p\!\left(\frac{y-a}{1-a}\right),
&y\in[a,2a);\cr
A_1(y+a)=A_1(y)+p(y+a)-\frac{1}{1-a}\;p\!\left(\frac{y}{1-a}\right),&y\in[a,1-2a);\cr
A_1(y)=\frac{1}{1-a}\;p\!\left(\frac{y}{1-a}\right)-p(y+a),&
y\in[1-2a,1-a].\cr
\end{cases}\end{equation}

Note that for $0\le y<a$ or $1-a<y\le 1$ function $\alpha_1(y)$ can
be arbitrary, because equality~\eqref{lemma1bis} doesn't apply any
conditions on it.

Let $\tilde{y}\in[a,2a)$, then
\begin{multline*}A_1(\tilde{y})=p(\tilde{y}-a)+p(\tilde{y})-\frac{1}{1-a}\;p\!\left(\frac{\tilde{y}-a}{1-a}\right)=\\=\sum\limits_{k=-1}^{0}p(\tilde{y}+ka)-\frac{1}{1-a}\sum\limits_{k=-1}^{-1}p\!\left(\frac{\tilde{y}+ka}{1-a}\right).\end{multline*}
Using the second equality of the system~\eqref{lemma2bis} we can get
by induction:
\begin{multline*}A_1(\tilde{y}+a)=\sum\limits_{k=-1}^{0}p(\tilde{y}+ka)-\frac{1}{1-a}\sum\limits_{k=-1}^{-1}p\!\left(\frac{\tilde{y}+ka}{1-a}\right)+p(\tilde{y}+a)-\frac{p\!\left(\frac{\tilde{y}}{1-a}\right)}{1-a}=\cr=\sum\limits_{k=-1}^{1}p(\tilde{y}+ka)-\frac{1}{1-a}\sum\limits_{k=-1}^{0}p\!\left(\frac{\tilde{y}+ka}{1-a}\right),\end{multline*}
\begin{multline*}A_1(\tilde{y}+2a)=\sum\limits_{k=-1}^{1}p(\tilde{y}+ka)-\frac{1}{1-a}\sum\limits_{k=-1}^{0}p\!\left(\frac{\tilde{y}+ka}{1-a}\right)+p(\tilde{y}+2a)-\frac{p\!\left(\frac{\tilde{y}+a}{1-a}\right)}{1-a}=\cr=\sum\limits_{k=-1}^{2}p(\tilde{y}+ka)-\frac{1}{1-a}\sum\limits_{k=-1}^{1}p\!\left(\frac{\tilde{y}+ka}{1-a}\right),\end{multline*}
$$\vdots$$
\begin{equation}\label{lemma4}A_1(\tilde{y}+ma)=\sum\limits_{k=-1}^{m}p(\tilde{y}+ka)-\frac{1}{1-a}\sum\limits_{k=-1}^{m-1}p\!\left(\frac{\tilde{y}+ka}{1-a}\right),\end{equation}
where $m=1,2,\ldots$ is such that $\tilde{y}+ma-a\in [a,1-2a)$, that
is
$$\tilde{y}\in [(2-m)a,1-(m+1)a)\cap [a,2a).$$ This way the values of function $A_1(y)$ from interval $[a,2a)$
are induced onto $[2a,3a)$, $[3a,4a)$,\ldots\ The
system~\eqref{lemma2bis} can be described by scheme depicted on
fig.~\ref{picex22}, in which $\gamma=1-\left[\frac{1-a}{a}\right]a$
($[x]$
--- whole part of $x\in\mathbb R$).
\begin{figure}[h]
\centering\includegraphics[width=\textwidth]{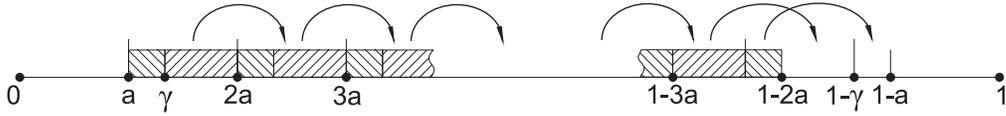}
\caption{Scheme of induction in
system~\eqref{lemma2bis}}\label{picex22}
\end{figure}

Let $\frac{1}{n+1}<a\le\frac{1}{n}$, $n=3,4,\ldots$ Note that
$\gamma=1-(n-1)a$ ($n\le\frac{1}{a}<n+1$,
$\left[\frac{1}{a}\right]=n$,
$\left[\frac{1-a}{a}\right]=\left[\frac{1}{a}\right]-1=n-1$) and
$a\le\gamma< 2a$. We will use system~\eqref{lemma2bis} again. In
order to improve visibility we will highlight interesting expression
in a chain of inequalities in bold. Take $y_0\in [a,\gamma)$, then,
for $m=n-2$,
$$a\le(n-2)a=a+(n-3)a\le\mathbf{y_0+(n-2)a-a}<1-(n-1)a+(n-3)a= 1-2a.$$

On the other hand, whereas
$$1-2a<(n-1)a\le\mathbf{y_0+(n-2)a}< 1-a,$$ we can use the formula~\eqref{lemma4} and the third equality of system~\eqref{lemma2bis}:
\begin{multline*}\frac{1}{1-a}\;p\!\left(\frac{y_0+(n-2)a}{1-a}\right)-p(y_0+(n-1)a)=A_1(y_0+(n-2)a)=\cr=\sum\limits_{k=-1}^{n-2}p(y_0+ka)-\frac{1}{1-a}\sum\limits_{k=-1}^{n-3}p\!\left(\frac{y_0+ka}{1-a}\right),\end{multline*}
which implies formula~\eqref{eqth1}.

Now take $y_1\in [\gamma,2a)$. If $n\ge4$, then, for $m=n-3$,
$$a\le1-3a=\gamma+(n-3)a-a\le\mathbf{y_1+(n-3)a-a}<2a+(n-3)a-a\le
1-2a.$$

On the other hand, whereas
$$1-2a\le\mathbf{y_1+(n-3)a}< 1-a,$$ we can write:
\begin{multline*}\frac{1}{1-a}\;p\!\left(\frac{y_1+(n-3)a}{1-a}\right)-p(y_1+(n-2)a)=A_1(y_1+(n-3)a)=\cr=\sum\limits_{k=-1}^{n-3}p(y_1+ka)-\frac{1}{1-a}\sum\limits_{k=-1}^{n-4}p\!\left(\frac{y_1+ka}{1-a}\right),\end{multline*}
which implies formula~\eqref{eqth2}.

If $n=3$, then $y_1\in [a,2a)\cap [1-2a,1-a]$. Therefore, we get
formula~\eqref{eqth2} by equating the first and the third equalities
in the system~\eqref{lemma2bis}.

 2. Now let  $\frac{1}{3}<a\le
\frac{1}{2}$. In this case $\frac{1-2a}{1-a}<\frac{a}{1-a}$ and
formula~\eqref{lemma1bis} can be written in the following way:
\begin{equation}\label{th2}\frac{p(x)}{1-a}=\begin{cases}p((1-a)x)+A_2((1-a)x+a),&
x\in\left[\left.0,\frac{1-2a}{1-a}\right.\right);\\[0.2cm]
p((1-a)x)+p((1-a)x+a),&x\in\left[\left.\frac{1-2a}{1-a},\frac{a}{1-a}\right.\right);\\[0.2cm]
p((1-a)x+a)+A_1((1-a)x),& x\in\left[\frac{a}{1-a},1\right].
\end{cases}\end{equation}

Consider the first equality of the system:
$\frac{p(x)}{1-a}=p((1-a)x)+A_2((1-a)x+a)$,
$x\in\left[\left.0,\frac{1-2a}{1-a}\right.\right)$. Make a change
$y=(1-a)x+a$, then, for $y\in[a,1-a)$,
\begin{equation}\label{newlabel}\frac{1}{1-a}\;p\!\left(\frac{y-a}{1-a}\right)=p(y-a)+A_2(y)=p(y-a)+p(y)-A_1(y).\end{equation}

Consider the second equality of the system:
$\frac{p(x)}{1-a}=p((1-a)x)+p((1-a)x+a)$,
$x\in\left[\left.\frac{1-2a}{1-a},\frac{a}{1-a}\right.\right)$. Make
a change $y=(1-a)x+a$, then, for $y\in[1-a,2a)$,
\begin{equation}\label{newlabel3}\frac{1}{1-a}\;p\!\left(\frac{y-a}{1-a}\right)=p(y)+p(y-a).\end{equation}

Consider the third equality of the system:
$\frac{p(x)}{1-a}=p((1-a)x+a)+A_1((1-a)x)$, $
x\in\left[\frac{a}{1-a},1\right]$. Make a change $y=(1-a)x$, then,
for $y\in[a,1-a]$,
\begin{equation}\label{newlabel2}A_1(y)=\frac{1}{1-a}\;p\!\left(\frac{y}{1-a}\right)-p(y+a).\end{equation}

Equating the values of $A_1(y)$ obtained from
equations~\eqref{newlabel} and~\eqref{newlabel2} for $y\in[a,1-a)$,
and taking in account equation~\eqref{newlabel3}, we exactly get
formulae~\eqref{eqth1}--\eqref{eqth3} for $n=2$.






Conversely, $\mu_S=\mu$ provided that~\eqref{eqth1}--\eqref{eqth3}
hold true. Indeed, these conditions do not set a value of
$\alpha_1(x)$ in the point $1-a$, but on the interval $[a,1-a)$ they
are equivalent to the system~\eqref{lemma2bis} (for  $n\ge 3$) or to
conditions~\eqref{newlabel}--\eqref{newlabel2} (for $n=2$).
Therefore almost everywhere (excluding point$x=1$ ($y=1-a$))
equality~\eqref{lemma1bis} holds true.
\end{proof}

For the sake of clarity we provide two Corollaries from the Theorem.

\begin{Sequence}Given measure $\mu\ll\lambda$, there exists equipped $2$-transformation $S(a,\alpha)$ preserving measure $\mu$ if and only if parameters $p(x)$,
$a$ and $\alpha=\{\alpha_1(x),\alpha_2(x)\}$ suffice the
conditions~\eqref{eqth1}--\eqref{eqth3}.
\end{Sequence}

\begin{Sequence}
Given equipped $2$-transformation $S(a,\alpha)$, there exists
measure $\mu\ll\lambda$ which is preserved by $S$ if and only if the
parameters $p(x)$, $a$ and $\alpha=\{\alpha_1(x),\alpha_2(x)\}$
suffice the conditions~\eqref{eqth1}--\eqref{eqth3}.
\end{Sequence}

An example of trivial density is given in the following Theorem
discussing the case of the Lebesgue measure ($\mu=\lambda$,
$p\equiv1$).

\begin{Theorem}\label{lebesguemeasure}
$\lambda_S=\lambda$ if and only if $a=\frac{1}{n}$, $n\in\mathbb N$,
$n\ge2$, and
$$\alpha_1(x)=\begin{cases}
\frac{n-2}{n-1},&x\in\left[\left.\frac{1}{n},\frac{2}{n}\right.\right);\cr
\frac{n-3}{n-1},&x\in\left[\left.\frac{2}{n},\frac{3}{n}\right.\right);\cr
\hspace{2.9 mm}\vdots\cr \frac{1}{n-1},&
x\in\left[\left.\frac{n-2}{n},\frac{n-1}{n}\right.\right).
\end{cases}$$
On the intervals $\left[0,\frac{1}{n}\right)$,
$\left[\frac{n-1}{n},1\right]$ function $\alpha_1(x)$ is arbitrary.
\end{Theorem}

\begin{proof}
We apply Theorem~\ref{main} for $p(x)=1$. The
condition~\eqref{eqth1} loses its meaning, because the interval
turns into an empty set:
$$n+1=\frac{1}{1-a}n,\quad x_0\in[a,1-(n-1)a)=\left[\frac{1}{n+1},\frac{1}{n+1}\right)=\emptyset.$$
Condition~\eqref{eqth2},
$$n=\frac{1}{1-a}(n-1),\quad
x_1\in[1-(n-1)a,2a)=\left[\frac{1}{n},\frac{2}{n}\right),$$
immediately implies $a=\frac{1}{n}$. Condition~\eqref{eqth3} for
$m=\overline{0,n-3}$ yields the following formula:
$$\alpha_1\!\left(x+\frac{m}{n}\right)=m+2-\frac{m+1}{1-\frac{1}{n}}=\frac{n-m-2}{n-1},\quad x\in\left[\frac{1}{n},\frac{2}{n}\right).$$

For $m=n-2$ the Condition~\eqref{eqth3} loses its meaning, because
$[(n-1)a,1-a)=\emptyset$. In particular for $n=2$ the equipment can
be chosen arbitrarily.

Thus summing up aforesaid we obtain the formula
$$\alpha_1(x)=\frac{n-k}{n-1},\quad x\in\left[\frac{k-1}{n},\frac{k}{n}\right),\quad k=\overline{2,n-1},$$
which implies the statement of the Theorem.
\end{proof}

\begin{Sequence}\label{diadic} If  $a=\frac{1}{2}$, then  $\lambda_S=\lambda$ for every $\alpha$.
\end{Sequence}

Corollary~\ref{diadic} agrees with a known result~\cite{renyi1957}:
diadic transformation $S(x)=2x\ (\!\!\!\mod 1)$ preserves the
Lebesgue measure.

 However, our construction admits not only trivial density: equalities~\eqref{eqth1}--\eqref{eqth3} are possible when $p(x)$ is not a constant.

\begin{Theorem}\label{interesting}
For all $n=2,3,\ldots$ there exists parameter $a$,
$\frac{1}{n+1}<a<\frac{1}{n}$, density  $p(x)$ and equipment
$\{\alpha_1(x),\alpha_2(x)\}$ such that $\mu_S=\mu$. Furthermore
density $p(x)$ is not a constant.
\end{Theorem}

\begin{proof}
First consider the case of even $n=2m$. Let
$a=\frac{2m+1-\sqrt{4m^2+1}}{2m}=\frac{n+1-\sqrt{n^2+1}}{n}$ be the
root of the equation
\begin{equation}\label{auxilliary}\frac{ma}{1-a}=1-ma.\end{equation}
It is easy to verify that
\begin{equation}\label{aux2}\frac{1}{n+1}<a<\frac{1}{n}\qquad \left(\frac{1}{2m+1}<a<\frac{1}{2m}\right).\end{equation}

Let $\beta,\gamma\ge0$. Consider a density given by the step
function
 $$p(x)=\beta\chi_{[0,ma)}(x)+(\beta+\gamma)(1-ma)\chi_{[ma,1-ma)}(x)+\gamma\chi_{[1-ma,1]}(x).$$
  Let us show that it suffices the equalities~\eqref{eqth1}--\eqref{eqth2}.
  Thereto henceforward we will use conditions~\eqref{auxilliary}--\eqref{aux2}, and for better visibility will also highlight interesting parts of equalities in bold.

  Consider condition~\eqref{eqth1}. Let $x_0\in[a,1-(2m-1)a)$. For
  $-1\le k\le m-2$
  $$0\le\mathbf{\frac{x_0+ka}{1-a}}<\frac{1-(2m-1)a+ka}{1-a}\le\frac{1-(m+1)a}{1-a}=ma.$$

   For $m-1\le k\le n-2$
  $$1-ma=\frac{ma}{1-a}\le\frac{a+ka}{1-a}\le\mathbf{\frac{x_0+ka}{1-a}}<1.$$

  For
  $-1\le k\le m-2$
  $$0\le\mathbf{ x_0+ka}<1-(2m-1)a+ka<(2m+1)a-(2m-1)a+ka=a(2+k)\le ma.$$

 For $k=m-1$
  $$ma=(k+1)a\le\mathbf{ x_0+ka}<1-(2m-1)a+ka= 1-ma.$$

 For $m\le k\le n-1$
  $$1-ma<(2m+1)a-ma=ma+a\le a+ak\le\mathbf{ x_0+ka}<1.$$

  In the case $\beta=\gamma=0$ equalities~\eqref{eqth1}--\eqref{eqth2} are obviously fulfilled. Henceforth, we will assume $\beta+\gamma>0$.
  Then
  $$\frac{\sum\limits_{k=-1}^{n-2}p\!\left(\frac{x_0+ka}{1-a}\right)}{\sum\limits_{k=-1}^{n-1}p(x_0+ka)}=\frac{m\beta+m\gamma}{m\beta+(\beta+\gamma)(1-ma)+m\gamma}=\frac{m}{1-ma+m}=1-a.$$

  Consider condition~\eqref{eqth2}. Let $x_1\in[1-(2m-1)a,2a)$. For
  $-1\le k\le m-3$
  $$0<\mathbf{\frac{x_1+ka}{1-a}}<\frac{2a+ka}{1-a}\le\frac{ma-a}{1-a}<ma.$$

For $k= m-2$
\begin{multline*}ma=\frac{1-(m+1)a}{1-a}=\frac{1-2ma+a+ka}{1-a}\le\mathbf{\frac{x_1+ka}{1-a}}<\frac{2a+ka}{1-a}=\frac{ma}{1-a}=\\=1-ma.\end{multline*}

 For $m-1\le k\le n-3$
  $$1-ma<\frac{1-2ma+a+ma-a}{1-a}\le\frac{1-(2m-1)a+ka}{1-a}\le\mathbf{\frac{x_1+ka}{1-a}}<1.$$

For $-1\le k\le m-2$
$$0<\mathbf{x_1+ka}<2a+ka\le ma.$$

For $m-1\le k\le n-2$
$$1-ma=1-(2m-1)a+ma-a\le\mathbf{ x_1+ka}<1.$$

Thus
 $$\frac{\sum\limits_{k=-1}^{n-3}p\!\left(\frac{x_1+ka}{1-a}\right)}{\sum\limits_{k=-1}^{n-2}p(x_1+ka)}\!=\!\frac{(m-1)\beta+(\beta+\gamma)(1-ma)+(m-1)\gamma}{m\beta+m\gamma}\!=\!\frac{m(1-a)}{m}\!=\!1-a.$$

The case of odd $n=2m-1$, $m=2,3,\ldots$ we will describe in less
detail. Let
$a=\frac{m-\sqrt{m^2-m}}{m}=\frac{n+1-\sqrt{n^2-1}}{n+1}$ be the
root of the equation
\begin{equation}\label{auxilliary2}\frac{(m-1)a}{1-a}=1-ma\qquad
\left(\frac{1}{n+1}<a<\frac{1}{n}\right).\end{equation}

Let
$p(x)=\beta\chi_{[0,1-ma)}(x)+(\beta+\gamma)(1-ma)\chi_{[1-ma,ma)}(x)+\gamma\chi_{[ma,1]}(x)$,
where $\beta,\gamma\ge0$.
  Let us show that $p(x)$ suffices the equalities~\eqref{eqth1}--\eqref{eqth2}.

  Consider condition~\eqref{eqth1}. Let $x_0\in[a,1-(2m-2)a)$.
 Then
  $$\hskip3.3mm\frac{x_0+ka}{1-a}\in\begin{cases}[0,1-ma),& -1\le k\le m-3;\cr
  [1-ma,ma),&k=m-2;\cr
  (ma,1),&m-1\le k\le n-2,\end{cases}$$
 $$x_0+ka\in\begin{cases}[0,1-ma),& -1\le k\le m-2;\cr
  [ma,1),&m-1\le k\le n-1,\end{cases}$$
$$\frac{\sum\limits_{k=-1}^{n-2}p\!\left(\frac{x_0+ka}{1-a}\right)}{\sum\limits_{k=-1}^{n-1}p(x_0+ka)}\!=\!\frac{(m-1)\beta+(\beta+\gamma)(1-ma)+(m-1)\gamma}{m\beta+m\gamma}\!=\!\frac{m(1-a)}{m}\!=\!1-a.$$



Consider condition~\eqref{eqth2}. Let $x_1\in[1-(2m-2)a,2a)$.
  Then
  $$\frac{x_1+ka}{1-a}\in\begin{cases}(0,1-ma),& -1\le k\le m-3;\cr
  [ma,1),&m-2\le k\le n-3,\end{cases}\phantom{\ \;\,}$$
 $$x_1+ka\in\begin{cases}(0,1-ma),& -1\le k\le m-3;\cr
 [1-ma,ma),& k=m-2;\cr
  (ma,1),&m-1\le k\le n-2,\end{cases}$$
  $$\frac{\sum\limits_{k=-1}^{n-3}p\!\left(\frac{x_1+ka}{1-a}\right)}{\sum\limits_{k=-1}^{n-2}p(x_1+ka)}=\frac{(m-1)\beta+(m-1)\gamma}{(m-1)\beta+(\beta+\gamma)(1-ma)+(m-1)\gamma}=\frac{m-1}{m-ma}=1-a.$$


It remains to show that for the function $\alpha_1\!(x)$ given by
the formula~\eqref{eqth3},  $0\!\le\!\alpha_1(x)\!\le\!1$. We will
find out $\alpha_1(x)$ explicitly.

Recall formula~\eqref{eqth3}:
\begin{equation*}\alpha_1(x+sa)p(x+sa)=\sum\limits_{k=-1}^{s}p(x+ka)-\frac{1}{1-a}\sum\limits_{k=-1}^{s-1}p\!\left(\frac{x+ka}{1-a}\right),\end{equation*}
where $x+sa\in\left\{
      \begin{array}{ll}
        [(s+1)a,(s+2)a) & \hbox{for } s=\overline{0,n-3}; \cr
        [(s+1)a,1-a) & \hbox{for } s=n-2.
      \end{array}
    \right.$

To find $\alpha_1(x)$ we have to consider $10$ following cases. In
doing this, we will use inequalities from above on numbers $x+ka$
and $\frac{x+ka}{1-a}$. In cases 1.1.1--2.2.3, when we divide by
$\beta$, $\gamma$ or $\beta+\gamma$, we assume $\beta>0$, $\gamma>0$
or $\beta+\gamma>0$ correspondingly. Otherwise ($\beta=0$,
$\gamma=0$ or $\beta+\gamma=0$) the values of $\alpha_1(x+sa)$ can
be chosen arbitrarily ($\alpha_1(x)\in L^1$, $0\le\alpha_1(x)\le1$).

1. $n=2m$. Note here that $\frac{ma}{1-a}=1-ma$.

1.1. $x\in[a,1-(n-1)a)$.

1.1.1. $0\le s\le m-2$, then
$$\alpha_1(x+sa)=\frac{1}{\beta}\left((s+2)\beta-\frac{1}{1-a}(s+1)\beta\right)=s+2-\frac{s+1}{1-a}.$$

1.1.2. $s=m-1$, then
\begin{multline}\label{particularinteresting}\alpha_1(x+sa)=\\=\frac{1}{(\beta+\gamma)(1-ma)}\left((s+1)\beta+(\beta+\gamma)(1-ma)-\frac{1}{1-a}(s+1)\beta\right)\!\!=\frac{\gamma}{\beta+\gamma}.\end{multline}

1.1.3. $m\le s\le n-2$, then
\begin{multline*}\alpha_1(x+sa)=\\=\frac{1}{\gamma}\left(m\beta+(\beta+\gamma)(1-ma)+(s-m+1)\gamma-\frac{1}{1-a}(m\beta+(s-m+1)\gamma)\right)\!\!=\\=\frac{1}{\gamma}\left((1-ma)\gamma+(s-m+1)\gamma-\frac{1}{1-a}(s-m+1)\gamma\right)=\\=\frac{a(2m-s-1)}{1-a}.\end{multline*}

1.2. $x\in[1-(n-1)a,2a)$.

1.2.1. $0\le s\le m-2$, then
$$\alpha_1(x+sa)=\frac{1}{\beta}\left((s+2)\beta-\frac{1}{1-a}(s+1)\beta\right)=s+2-\frac{s+1}{1-a}.$$


1.2.2. $m-1\le s\le n-3$, then
\begin{multline*}\alpha_1(x+sa)=\frac{1}{\gamma}\left(\hskip-0.9cm\phantom{\frac{1}{1-a}}m\beta+(s-m+2)\gamma-\right.\\\left.-\frac{1}{1-a}((m-1)\beta+(\beta+\gamma)(1-ma)+(s-m+1)\gamma)\right)=\\=\frac{a(2m-s-2)}{1-a}.\end{multline*}

2. $n=2m-1$. Note here that $\frac{(m-1)a}{1-a}=1-ma$.

2.1. $x\in[a,1-(2m-2)a)$.

2.1.1. $0\le s\le m-2$, then
$$\alpha_1(x+sa)=\frac{1}{\beta}\left((s+2)\beta-\frac{1}{1-a}(s+1)\beta\right)=s+2-\frac{s+1}{1-a}.$$

2.1.2. $m-1\le s\le n-2$, then
\begin{multline*}\alpha_1(x+sa)=\frac{1}{\gamma}\left(m\beta+(s-m+2)\gamma-\phantom{\frac{1}{1-a}}\right.\\\left.-\frac{1}{1-a}((m-1)\beta+(\beta+\gamma)(1-ma)+(s-m+1)\gamma\right)=\frac{a(2m-s-2)}{1-a}.\end{multline*}

2.2. $x\in[1-(2m-2)a,2a)$.

2.2.1. $0\le s\le m-3$, then
$$\alpha_1(x+sa)=\frac{1}{\beta}\left((s+2)\beta-\frac{1}{1-a}(s+1)\beta\right)=s+2+\frac{s+1}{1-a}.$$

2.2.2. $s=m-2$, then
\begin{multline*}\alpha_1(x+sa)=\\=\frac{1}{(\beta+\gamma)(1-ma)}\left((s+1)\beta+(\beta+\gamma)(1-ma)-\frac{1}{1-a}(s+1)\beta\right)=\frac{\gamma}{\beta+\gamma}.\end{multline*}

2.2.3. $m-1\le s\le n-3$, then
\begin{multline*}\alpha_1(x+sa)=\frac{1}{\gamma}\left((m-1)\beta+(\beta+\gamma)(1-ma)+(s-m+2)\gamma-\phantom{\frac{1}{1-a}}\right.\\\left.-\frac{1}{1-a}((m-1)\beta+(s-m+2)\gamma)\right)=\frac{a(2m-s-3)}{1-a}.\end{multline*}

We sum up formulae for $\alpha_1(x)$ obtained above in more compact
way.

For $n=2m$, $x\in[a,1-(n-1)a)$,
$$\alpha_1(x+sa)=\left\{
                   \begin{array}{ll}
                     s+2-\frac{s+1}{1-a}, & 0\le s\le m-2; \\
                     \frac{\gamma}{\beta+\gamma}, & s=m-1; \\
                     \frac{a(n-s-1)}{1-a}, & m\le s\le
                     n-2.\phantom{\hspace{0.5cm}}
                   \end{array}
                 \right.
$$

For $n=2m$, $x\in[1-(n-1)a,2a)$,
$$\alpha_1(x+sa)=\left\{
                   \begin{array}{ll}
                     s+2-\frac{s+1}{1-a}, & 0\le s\le m-2; \\
                     \frac{a(n-s-2)}{1-a}, & m-1\le s\le n-3.
                   \end{array}
                 \right.
$$

For $n=2m-1$, $x\in[a,1-(n-1)a)$,
$$\alpha_1(x+sa)=\left\{
                   \begin{array}{ll}
                     s+2-\frac{s+1}{1-a}, & 0\le s\le m-2; \\
                     \frac{a(n-s-1)}{1-a}, & m-1\le s\le n-2.
                   \end{array}
                 \right.
$$

For $n=2m-1$, $x\in[1-(n-1)a,2a)$,
$$\alpha_1(x+sa)=\left\{
                   \begin{array}{ll}
                     s+2-\frac{s+1}{1-a}, & 0\le s\le m-3; \\
    \frac{\gamma}{\beta+\gamma}, & s=m-2; \\
                     \frac{a(n-s-2)}{1-a}, & m-1\le s\le n-3.
                   \end{array}
                 \right.
$$
Thus $0<\alpha_1(x)<1$, since for $k=1,2$
$$0<\frac{1-ma}{1-a}\le\mathbf{s+2-\frac{s+1}{1-a}}=\frac{1-2a-sa}{1-a}\le \frac{1-2a}{1-a}<1,\quad 0\le s\le m-2,$$
$$0<\frac{a}{1-a}\le\mathbf{\frac{a(n-s-k)}{1-a}}\le\frac{am}{1-a}<1,\quad m-1\le s\le n-3.$$


Using Theorem~\ref{main} for an equipment $\alpha$ given by the
formula~\eqref{eqth3} we obtain $\mu_S=\mu$.
\end{proof}

\begin{Remark*}In the proof of Theorem~\ref{interesting}  we  have actually
found a family of densities depending on two parameters
$\beta,\gamma\ge0$.
\end{Remark*}

\begin{figure}[!htb]
\centering\includegraphics[width=6cm]{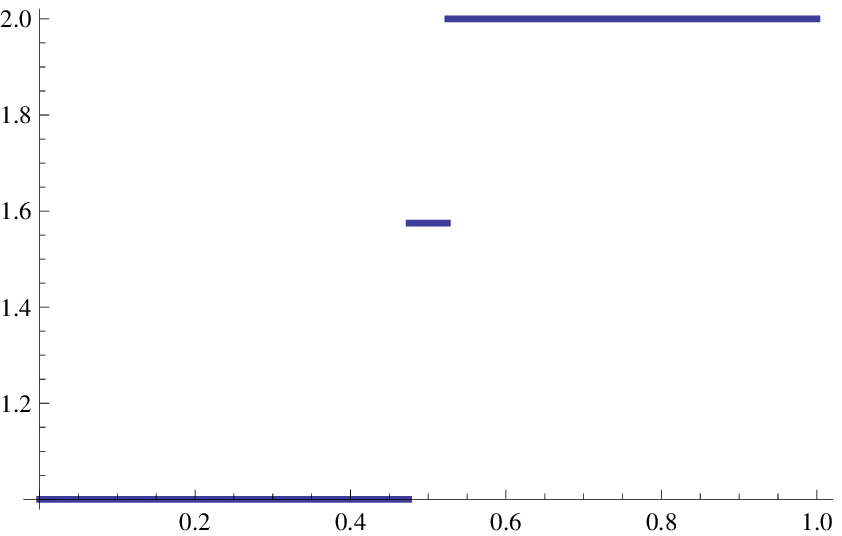}\quad\includegraphics[width=6cm]{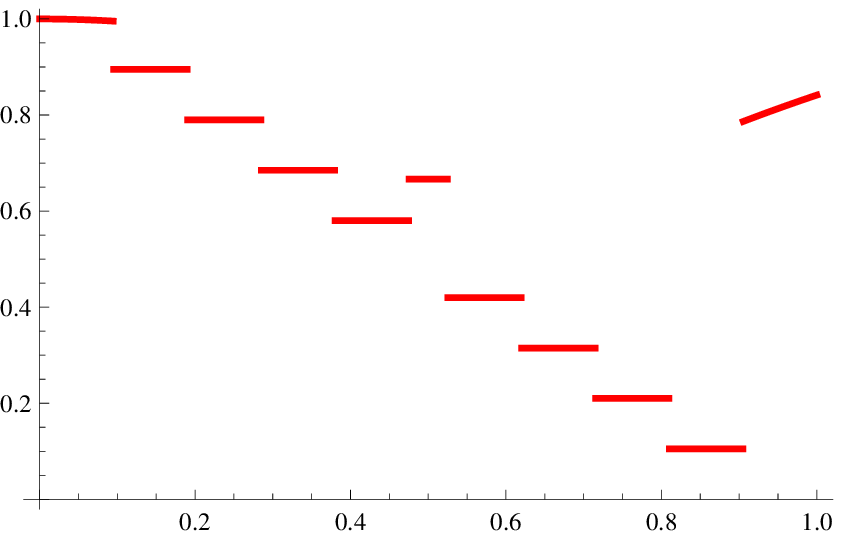}
\caption{Density  $p(x)$ plot (to the left) and weight function
$\alpha_1(x)$ plot (to the right), here $n=10$}\label{ex10}
\end{figure}

As an example we provide the plots of density  $p(x)$ and weight
function $\alpha_1(x)$ for the case of $n=10$, see fig.~\ref{ex10}.
We used $\beta=1$, $\gamma=2$, $\alpha_1(x)=\cos x$ for $x\in[0,a)$,
$\alpha_1(x)=\sin x$ for $x\in[1-a,1]$.

\begin{Sequence}
For all $\beta,\gamma\ge0$ chosen in the proof of
Theorem~\ref{interesting},
$$\mu([0,1])=\left\{
                                                    \begin{array}{ll}
                                                      \left(1+n^2-n\sqrt{n^2+1}\right)(\beta+\gamma), & n\hbox{ is even;} \\
                                                      \left(1-n^2+n\sqrt{n^2-1}\right)(\beta+\gamma), & n\hbox{ is odd.}
                                                    \end{array}
                                                  \right.
$$
In particular we can chose $\beta+\gamma$ such that $\mu([0,1])=1$.
\end{Sequence}

\begin{proof}
If $n=2m$, then $ma=(n+1-\sqrt{n^2+1})/2$ (see
formula~\eqref{auxilliary}) and
\begin{multline*}\!\!\!\mu([0,1])=\int_0^1\!\left(\beta\chi_{[0,ma)}(x)+
(\beta+\gamma)(1-ma)\chi_{[ma,1-ma)}(x)+\gamma\chi_{[1-ma,1]}(x)\right)d\lambda=\\=
\left(ma+(1-ma)(1-2ma)\right)(\beta+\gamma)=
\left(1+n^2-n\sqrt{n^2+1}\right)(\beta+\gamma).\end{multline*} If
$n=2m-1$, then $ma=(n+1-\sqrt{n^2-1})/2$ (see
formula~\eqref{auxilliary2}) and
\begin{multline*}\!\!\!\mu([0,1])=\int_0^1\!\left(\beta\chi_{[0,1-ma)}(x)+(\beta+\gamma)(1-ma)\chi_{[1-ma,ma)}(x)+\gamma\chi_{[ma,1]}(x)\right)d\lambda=\\=
\left(2ma(1-ma)\right)(\beta+\gamma)=
\left(1-n^2+n\sqrt{n^2-1}\right)(\beta+\gamma).\end{multline*}
\end{proof}

As a conclusion let us consider the following. Let $n=2$. By
Theorem~\ref{interesting} for the number
$a=\frac{3-\sqrt{5}}{2}\in\left(\frac{1}{3},\frac{1}{2}\right]$ and
density $p(x)=\beta\chi_{[0,a)}(x)+\beta(1-a)\chi_{[a,1-a)}(x)$ (we
use $\gamma=0$)
 $\mu_S=\mu$. Condition~\eqref{eqth3} turns to the following:
$$\alpha_1(x)p(x)=p(x-a)+p(x)-\frac{1}{1-a}\;p\!\left(\frac{x-a}{1-a}\right),\quad  x\in[a,1-a).$$

Then $x-a\in[0,1-2a)\subset[0,a)$,
$\frac{x-a}{1-a}\in\left[0,\frac{1-2a}{1-a}\right)=[0,a)$ and we
obtain
$$\alpha_1(x)\beta=\frac{1}{1-a}\left(\beta+\beta(1-a)-\frac{1}{1-a}\beta\right)=0$$
(this result also follows from formula~\eqref{particularinteresting}
given in the proof of Theorem~\ref{interesting}).

\begin{figure}[h]
\centering\includegraphics[clip,height=3.5cm]{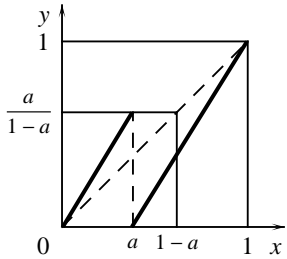}
\caption{Scheme of transformations $S_2$ and $Tx=\frac{1}{1-a}x
(\!\!\!\mod 1-a)$\hspace{3cm} (in the small square
region)}\label{ergodictransformation}
\end{figure}

Let us take $\alpha_1(x)\equiv0$, $\alpha_2(x)\equiv1$. This case
corresponds to single-valued dynamical system $([0,1],\mathfrak B,
\mu,S_2)$ consisting of one function $S_2$. The support of a measure
of such system lays inside $[0,1-a]$. Being restricted to $[0,1-a]$,
transformation $S_2$ coincides with the mapping
$$Tx=\frac{1}{1-a}x \quad(\!\!\!\!\!\mod 1-a),$$ acting by the rule $Tx=(1-a)\left\{\frac{1}{(1-a)^2}x\right\}$.
Note that $\frac{1}{1-a}=\frac{\sqrt{5}+1}{2}$ is a <<golden ratio>>
(being also a  <<simple $\beta$-number>>, see~\cite{parry1960}). On
the fig.~\ref{ergodictransformation} we depict  the plot of
transformation $S_2$ in the square $[0,1]\times[0,1]$ and the plot
of mapping $T$ in the smaller square $[0,1-a]\times[0,1-a]$.

Let us normalize the density $p(x)$ such that $\mu([0,1-a])=1-a$:
$$\mu([0,1-a])=\int_{0}^{1-a}p(x)\,d\lambda=\beta a+\beta (1-a)(1-2a)=1-a,$$
whence $\beta=\frac{5+3\sqrt{5}}{10}$,
$\beta(1-a)=\frac{5+\sqrt{5}}{10}$. For the sake of simplicity we
let also $p(1-a)=\frac{5+\sqrt{5}}{10}$ (since p$(x)$ is defined
$\lambda$-almost everywhere).

Thus we get invariant measure for the dynamical system
$([0,1-a],\mathfrak B',\mu', T)$ (here $\mathfrak B'$ is a
$\sigma$-field obtained by intersecting  sets from $\mathfrak B$
with the segment $[0,1-a]$, $\mu'=\mu_{|_{\mathfrak B'}}$).

\begin{Sequence}\label{invmeasurecorollary} Transformation $S_2|_{[0,1-a]}(x)=Tx\colon
[0,1-a]\to[0,1-a]$ preserves measure $\mu'(\tilde{p})$ with density
$$\tilde{p}(x)=\left\{
         \begin{array}{ll}
           \frac{5+3\sqrt{5}}{10}, & 0\le x<\frac{3-\sqrt{5}}{2}; \\
           \frac{5+\sqrt{5}}{10}, & \frac{3-\sqrt{5}}{2}\le x\le \frac{\sqrt{5}-1}{2} \hbox{.}
         \end{array}
       \right.$$
\end{Sequence}

This agrees with classical result (adapted for the segment
$[0,1-a]$) obtained by A.~R\'enyi in~\cite{renyi1957} (see
also~\cite{parry1960}): transformation $T$ is ergodic and has a
unique invariant measure (with density $\tilde{p}(x)$) equivalent to
the  Lebesgue measure.


\References{
\bibitem{renyi1957}
{\it A.~R\'enyi}, ``Representations for real numbers and their
ergodic properties,''
  {\em Acta Math. Acad. Sci. Hung.}, vol.~8, pp.~477--493, 1957.
\bibitem{parry1960}
{\it W.~Parry}, ``On the $\beta$-expansions of real numbers,'' {\em
Acta Math. Acad.
  Sci. Hung.}, vol.~11, pp.~401--416, 1960.
\bibitem{sidorov2003b}
{\it N.~Sidorov}, ``Arithmetic dynamics,'' in {\em Topics in
dynamics and ergodic
  theory} (S.~Bezuglyi and S.~F. Kolyada, eds.), vol.~310 of {\em London
  Mathematical Society lecture note series}, pp.~145--189, Cambridge University
  Press, 2003.
\bibitem{erdosh1990}
{\it P.~Erd\H{o}s, I.~Jo\'o, and V.~Komornik}, ``Characterization of
the unique
  expansions $1=\sum_{i=1}^\infty q^{-n_i}$ and related problems,'' {\em Bull.
  Soc. Math. France}, vol.~118, pp.~377--390, 1990.
\bibitem{sidorov2003}
{\it N.~Sidorov}, ``Almost every number has a continuum of
$\beta$-expansions,'' {\em
  Amer. Math. Monthly}, vol.~110, pp.~838--842, 2003.
\bibitem{sidorov2009}
{\it N.~Sidorov}, ``Expansions in non-integer bases: Lower, middle
and top orders,''
  {\em J. Number Theory}, vol.~129, no.~4, pp.~741--754, 2009.
\bibitem{KB2009}
{\it К.~Б. Игудесман}, ``Верхние адреса для одного семейства систем
итерированных
  функций на отрезке,'' {\em Известия вузов. Математика}, vol.~9, pp.~75--81,
  2009. (in Russian, but also available in English in {\em
Russian Mathematics (Iz. VUZ)}, vol.~53, no.~9, pp.~67--72, 2009.)
\bibitem{gelfond1959}
{\it А.~О. Гельфонд}, ``Об одном общем свойстве систем счисления,''
{\em Известия
  академии наук СССР. Серия математическая}, vol.~23, pp.~809--814, 1959.
\bibitem{barnsley2005}
{\it M.~F. Barnsley}, ``Theory and application of fractal tops,'' in
{\em Fractals in
  engineering: new trends in theory and applications} (J.~L\'evy-V\'ehel and
  E.~Lutton, eds.), pp.~3--20, Springer-Verlag, London Limited, 2005.
\bibitem{erdosh1939}
{\it P.~Erd\H{o}s}, ``On a family of symmetric Bernoulli
convolutions,'' {\em Amer.
  J. Math.}, vol.~61, pp.~974--975, 1939.
\bibitem{convolutions2000}
{\it Y.~Peres, W.~Schlag, and B.~Solomyak}, ``Sixty years of
Bernoulli
  convolutions,'' in {\em Fractal Geometry and Stochastics II} (C.~Bandt,
  S.~Graf, and M.~Zaehle, eds.), vol.~46 of {\em Progress in probability},
  pp.~39--65, Birkhauser, 2000.
\bibitem{troshin2009}
{\it P.~I. Troshin}, ``Multivalued dynamic systems with weights,''
{\em Izvestiya Vysshikh Uchebnykh Zavedenii. Matematika}, vol.~7,
pp.~35--50, 2009. (in Russian, but also available in English in {\em
Russian Mathematics (Iz. VUZ)}, vol.~53, no.~7, pp.~28--42, 2009.)}


 {\small \vspace{\baselineskip}\hrule \vspace{3pt}
\par
{\bf Paul I. Troshin} -- chair of Geometry of Kazan State
University, Kazan, Russia
\par
E-mail: {\it Paul.Troshin@gmail.com} }
\end{document}